\documentclass[a4paper,12pt]{article}
\usepackage{amsfonts}

\begin{document}
\newtheorem{theor}{Theorem}[section] 
\newtheorem{prop}[theor]{Proposition} 
\newtheorem{cor}[theor]{Corollary}
\newtheorem{lemma}[theor]{Lemma}
\newtheorem{sublem}[theor]{Sublemma}
\newtheorem{defin}[theor]{Definition}
\newtheorem{conj}[theor]{Conjecture}

\hfuzz2cm

\def\deg{{\widehat {\rm deg}\,}}
\def\bz{\mbox{\boldmath$\zeta$\unboldmath}}
\def\bzs{\mbox{\boldmath$\zeta'$\unboldmath}}
\def\bzo{\overline{\bz}}
\def\odd{{\rm odd}}
\def\a{\alpha}
\def\bC{{\bf C}}
\def\Td{{\rm Td}}
\def\ch{{\rm ch}}
\def\Hom{{\rm Hom}}
\def\Ad{{\rm Ad}^{1,0}_{G/P}}
\def\KO{H}
\def\K1{K}
\def\covol{{\rm covol}\,}
\gdef\beginProof{\par{\bf Proof: }}
\gdef\endProof{${\bf Q.E.D.}$\par}
\gdef\mtr#1{\overline{#1}}
\gdef\ar#1{\widehat{#1}}

\def \x{\times}
\def \BN{{\bf N}}
\def \BZ{{\bf Z}}
\def \BQ{{\bf Q}}
\def \BR{{\bf R}}
\def \BC{{\bf C}}
\def \BG{{\bf G}}
\def \<{\langle}
\def \>{\rangle}
\def \Spec {\hbox{\rm Spec }}
\def \X {{\rm Bor}(G)}
\def \Y {{\rm Par}(G)}
\def \ad {\hbox {ad}}

\author{
Kai K\"ohler\footnote{Centre de Math\'ematiques de Jussieu/C.P. 7012/2, place
Jussieu/F-75251 Paris Cedex 05/France/e-mail : koehler@math.jussieu.fr/URL:
http://www.math.jussieu.fr/$\tilde{\ }$koehler}} 
\title{A Hirzebruch proportionality principle in Arakelov geometry}
\maketitle
\begin{abstract}
We show that a conjectural extension of a fixed point formula in Arakelov geometry
implies results about a tautological subring in the arithmetic Chow ring of
bases of abelian schemes. Among the results are an Arakelov version of the Hirzebruch
proportionality principle and a formula for a critical power of $\ar c_1$ of the
Hodge bundle.
\end{abstract}
\begin{center}
\end{center}
\thispagestyle{empty}
\setcounter{page}{1}
\section{Introduction}
The purpose of this note is to exploit some implications of a conjectural fixed point
formula in Arakelov geometry when applied to the action of the $(-1)$ involution on
abelian schemes of relative dimension $d$. It it shown that the fixed point formula's
statement in this case is equivalent to giving the values of arithmetic Pontrjagin
classes of the Hodge bundle, where these Pontrjagin classes are defined as
polynomials in the arithmetic Chern classes defined by Gillet and Soul\'e. When
combined with the statement of the non-equivariant arithmetic
Grothendieck-Riemann-Roch formula, one obtains a formula for
the class $\hat c_1^{1+d(d-1)/2}$  of the Hodge bundle in terms of topological
classes and a certain special differential form $\gamma$. Finally we derive an
Arakelov version of the Hirzebruch proportionality principle, namely a ring
homomorphism from the Arakelov Chow ring of Lagrangian Grassmannians to the
arithmetic Chow ring of bases of abelian schemes.

A fixed point formula for maps from arithmetic varieties to Spec $D$
has been proven by Roessler and the author in \cite{KR1}, where $D$ is a regular
arithmetic ring. In
\cite[Appendix]{KR2} we described the generalization to flat equivariantly projective
maps between arithmetic varieties over $D$. The still missing ingredient to the proof
of this conjecture is the equivariant version of Bismut's formula for the behavior
of analytic torsion forms under the composition of immersions and fibrations
\cite{Bfam}, i.e. a merge of \cite{Bimm} and \cite{Bfam}.

We work only with regular schemes as bases; extending these results to moduli
stacks and their compactifications remains an open problem, as Arakelov geometry for
such situations is not yet developed. In particular one could search an analogue of
the full Hirzebruch-Mumford proportionality principle in Arakelov geometry. When this
article was almost finished, we learned about related work by van der Geer
concerning the classical Chow ring of the moduli stack of abelian varieties and its
compactifications
\cite{Geer}. The approach there to determine the tautological subring
uses the non-equivariant Grothendieck-Riemann-Roch theorem applied to the line
bundle associated to theta divisor. Thus it might be possible to avoid the use of the
fixed point formula in our situation by mimicking this method, possibly by extending
the methods of Yoshikawa \cite{Y}; but computing the occurring objects related to the
theta divisor is presumably not easy.

According to a conjecture by Oort, there are no complete subvarieties of codimension
$d$ in the complex moduli space for $d\geq3$. Thus a possible application of our
formula for $\hat c_1^{1+d(d-1)/2}$  of the Hodge bundle could be a proof of this
conjecture by showing that the height of potential subvarieties would be lower than
the known lower bounds for heights. Van der Geer \cite[Cor. 7.2]{Geer} used the
degree with respect to the Hodge bundle to show that complete subvarieties have
codimension
$\geq d$.

{\bf Acknowledgements: } I thank A. Johan de Jong, Damian Roessler, Christophe
Soul\'e and Thorsten Wedhorn for helpful discussions. Also I thank the Deutsche
Forschungsgemeinschaft which supported me with a Heisenberg fellowship during the
preparation of this article.

\section{Torsion forms}

Let $\pi:E^{1,0}\to B$ denote a $d$-dimensional holomorphic vector bundle over a
complex manifold. Let $\Lambda$ be a lattice subbundle of the  underlying
real vector bundle $E^{1,0}_\BR$ of rank
$2d$. Thus the quotient bundle $M:=E^{1,0}/\Lambda\to B$ is a holomorphic
fibration by tori $Z$. Let 
$$
\Lambda^*:=\{\mu\in (E^{1,0}_{\BR})^*\,|\,\mu(\lambda)\in2\pi{\bf
Z}\,\forall\lambda\in
\Lambda\}
$$
denote the dual lattice bundle. Assume that $E^{1,0}$ is equipped with an Hermitian
metric such that the volume of the fibres is constant. Such a metric is induced by a
polarization.

Let $N_V$ be the number operator acting on $\Gamma(Z,\Lambda^q
T^{*0,1}Z)$ by multiplication with $q$. Let ${\rm Tr}_s$ denote the supertrace with
respect to the $\BZ/2\BZ$-grading on $\Lambda T^*B\otimes{\rm End}(\Lambda
T^{*0,1}Z)$. Let
$\phi$ denote the map acting on
$\Lambda^{2 p}T^*B$ as multiplication by
$(2\pi i)^{-p}$.  We write 
$\widetilde{\cal A}(B)$ for 
$\widetilde{\cal A}(B):=\bigoplus_{p\geq 0}({\cal A}^{p,p}(B)/({\rm Im}\,\partial
+{\rm Im}\,\mtr{\partial}))$, where ${\cal A}^{p,p}(B)$ denotes the
$C^\infty$ differential forms of type $(p,p)$ on $B$.

In \cite[Section 3]{Ktori}, a superconnection $A_t$
acting on the infinite-dimensional vector bundle
$\Gamma(Z,\Lambda T^{*0,1}Z)$ over $B$ had been introduced, depending on $t\in\BR^+$.
For a fibrewise acting holomorphic isometry $g$ the limit $\lim_{t\to\infty} \phi
{\rm Tr}_s g^*N_H e^{-A_t^2}=:\omega_\infty$ exists and is given by the respective
trace restricted to the cohomology of the fibres. The equivariant analytic torsion
form
$T_g(\pi,{\cal O}_M)\in\widetilde{\cal A}(B)$ was defined there as the
derivative at zero of the zeta function with values in differential forms on $B$
given by
$$
-\frac1{\Gamma(s)}\int_0^\infty (\phi {\rm Tr}_s g^*N_H
e^{-A_t^2}-\omega_\infty)t^{s-1}\,dt
$$
for Re $s>d$.

\begin{theor} Let an isometry $g$ act fibrewise with isolated fixed points
on the fibration by tori $\pi:M\to B$. Then the equivariant torsion form
$T_g(\pi,{\cal O}_M)$ vanishes.
\end{theor}
\beginProof Let $f_\mu:M\to \BC$ denote the function $e^{i\mu}$ for
$\mu\in\Lambda^*$. As is shown in \cite[\S 5]{Ktori} the operator $A_t^2$ acts
diagonally with respect to the Hilbert space decomposition
$$
\Gamma(Z,\Lambda T^{*0,1}Z)=\bigoplus_{\mu\in\Lambda^*} \Lambda
E^{*0,1}\otimes\{f_\mu\}\,\,.
$$
As in \cite[Lemma 4.1]{KR4} the induced action by $g$
maps a function $f_\mu$ to a multiple of itself if only if $\mu=0$ because $g$ acts
fixed point free on $E^{1,0}$ outside the zero section. In that case, $f_\mu$
represent an element in the cohomology. Thus the zeta function defining the torsion
vanishes.
\endProof

{\bf Remark: }As in \cite[Lemma 4.1]{KR4}, the same proof shows the vanishing of the
equivariant torsion form $T_g(\pi,L)$ for coefficients in a $g$-equivariant line
bundle
$L$ with vanishing first Chern class.\medskip

We shall also need the following result of \cite{Ktori} for the non-equivariant
torsion form $T(\pi,{\cal O}_M):=T_{\rm id}(\pi,{\cal O}_M)$: Assume for simplicity
that
$\pi$ is K\"ahler. Consider for Re $s<0$ the zeta function with values in
$(d-1,d-1)$-forms on $B$
$$
Z(s):=\frac{\Gamma(2n-s-1)}{\Gamma(s)(n-1)!}
\sum_{\lambda\in\Lambda\setminus\{0\}}
(\frac{\mtr\partial\partial}{4\pi i}\|\lambda^{1,0}\|^2)^{\wedge(d-1)}
(\|\lambda^{1,0}\|^2)^{s+1-2n}
$$
where $\lambda^{1,0}$ denotes a lattice section in $E^{1,0}$. Then the limit
$\gamma:=\lim_{s\to 0^-}Z'(0)$ exists and
$$
\frac{\mtr\partial\partial}{2\pi i}\gamma=c_n(\mtr{E^{0,1}})\,\,.
$$
In \cite[Th. 4.1]{Ktori} the torsion form is shown to equal
$$
T(\pi,{\cal O}_M)=\frac\gamma{{\rm Td}(\mtr{E^{0,1}})}
$$
in $\widetilde{\cal A}(B)$. The differential form $\gamma$ was intensively studied in
\cite{Ktori}, in particular its behavior under the action of Hecke operators was
determined.

We shall denote a vector bundle $F$ together with an Hermitian metric $h$ by $\mtr
F$. Then
$\ch_g(\mtr{F})$ shall denote the Chern-Weil representative of the equivariant Chern
character
associated to the 
restriction of $(F,h)$ to the fixed point subvariety. Recall 
also that $\Td_{g}(\mtr{F})$ is the differential form 
$$\frac{\Td(\mtr{F}^g)}{\sum_{i\geq 0}(-1)^k
\ch_{g}(\Lambda^k\mtr{F})}\,\,.$$

\section{Abelian schemes and the fixed point formula}

We shall use the Arakelov geometric concepts and notation of \cite{Soule} and
\cite{KR1}. In this article we shall only give a brief introduction to Arakelov
geometry, and we refer to \cite{Soule} for details. Let
$D$ be a regular arithmetic ring, i.e. a regular,  excellent, Noetherian
integral ring, together  with a finite set $\cal S$ of ring monomorphism
of
$D\to{\bf C}$,  invariant under complex conjugation.
We shall denote by $\mu_n$ the diagonalizable group scheme over $D$ 
associated to ${\bf Z}/{n{\bf Z}}$. 
We choose once and for all 
a primitive $n$-th root of unity $\zeta_{n}$. 
Let $f:Y\to{\rm Spec}\ D$ be an equivariant arithmetic 
variety, i.e.
a regular integral scheme,  
endowed with a $\mu_n$-projective 
action over ${\rm Spec}\ D$.
The groups of $n$-th roots of unity acts 
on $Y(\BC)$ by holomorphic automorphisms and we shall 
write $g$ for the automorphism corresponding to $\zeta_{n}$.

We write $f^{\mu_n}$ for the
map $Y_{\mu_n}\to{\rm Spec}\ D$  induced by $f$ on the fixed point subvariety. 
Complex conjugation induces an antiholomorphic 
automorphism of $Y(\BC)$ and $Y_{\mu_{n},\BC}$, both of 
which we denote by $F_{\infty}$.  
The space 
$\widetilde{\cal A}(Y)$ is the subspace of $\widetilde{\cal A}(Y(\BC))$ of
classes of differential forms
$\omega$ such that $F_{\infty}^{*}\omega=(-1)^{p}\omega$.

Let $\ar{\rm CH}^*(Y)$ denote the Gillet-Soul\'e arithmetic Chow ring, consisting of
arithmetic cycles and suitable Green currents on $Y(\BC)$. Let ${\rm CH}^*(Y)$
denote the classical Chow ring. Then there is an exact sequence in any degree $p$
\begin{equation}
{\rm CH}^{p,p-1}(Y)\stackrel{\rho}{\to}
\widetilde{\cal A}^{p-1,p-1}(Y)\stackrel{a}{\to}
\ar{\rm CH}^p(Y)\stackrel{\zeta}{\to}
{\rm CH}^p(Y)\to0\,\,.
\end{equation}
For Hermitian vector bundles $\mtr E$ on $Y$ Gillet and Soul\'e defined arithmetic
Chern classes $\ar c_p(\mtr E)\in\ar{\rm CH}^*(Y)_\BQ$.

By "product of Chern classes", we shall understand in this article any product of at
least two equal or non-equal Chern classes of degree larger than 0 of a given vector
bundle.
\begin{lemma} \label{hups}
Let
$$
\ar\phi=\sum_{j=0}^\infty a_j \ar c_j+\mbox{\rm products of Chern classes}
$$
denote an arithmetic characteristic class with $a_j\in\BQ$, $a_j\neq0$ for $j>0$.
Assume that for a vector bundle
$\mtr F$ on an arithmetic variety $Y$, $\ar\phi(\mtr F)=m+a(\beta)$ where
$\beta$ is a differential form on $Y(\BC)$ with
$\partial\mtr\partial \beta=0$ and $m\in\ar{CH}^0(Y)_\BQ$. Then
$$
\sum_{j=0}^\infty a_j \ar c_j(\mtr F)=m+a(\beta)\,\,.
$$
\end{lemma}
\beginProof By induction: for the term in $\ar{CH}^0(Y)_\BQ$, the formula is clear.
Assume now for $k\in\BN_0$ that
$$
\sum_{j=0}^k a_j \ar c_j(\mtr F)=m+\sum_{j=0}^k a(\beta)^{[j]}\,\,.
$$
Then $\ar c_j(\mtr F)\in a({\rm ker}\,\partial\mtr\partial)$ for $1\leq j\leq
k$, thus products of these $\ar c_j$'s vanish by 
\cite[Remark III.2.3.1]{Soule}. Thus the term of degree $k+1$ of $\ar\phi(\mtr F)$
equals $a_{k+1} \ar c_{k+1}(\mtr F)$.
\endProof
We define {\bf arithmetic Pontrjagin classes} $\ar p_j\in \ar{CH}^{2j}$ of
arithmetic vector bundles by the relation
$$
\sum_{j=0}^\infty (-z^2)^j\ar p_j:=(\sum_{j=0}^\infty z^j\ar
c_j)(\sum_{j=0}^\infty (-z)^j\ar c_j)\,\,.
$$
Thus, $$
\ar p_j(\mtr F)=(-1)^j \ar c_{2j}(\mtr
F\oplus \mtr F^*)
=\ar c_j^2(\mtr F)+2\sum_{l=0}^{j-1}(-1)^{j+l}\ar c_l(\mtr F)\ar c_{2k-l}(\mtr F)
$$ for an arithmetic vector bundle $\mtr F$ (compare \cite[\S 15]{Milnor}).
Lemma \ref{hups} holds
with Chern classes replaced by Pontrjagin classes.

Let $f:Y\to\Spec D$ denote a quasi-projective arithmetic variety and let $\pi:X\to Y$
denote a principally polarized  abelian scheme of relative dimension $d$. Set $\mtr
E:=(R^1\pi_* {\cal O},\|\cdot\|_{L^2})^*$.
This bundle
$E={\bf Lie}(X/Y)^*$ is the Hodge bundle.
 Then by \cite[Prop. 2.5.2]{BBM}, the full
direct image of
$\cal O$ under $\pi$ is given by
\begin{equation}
\mtr{R^\bullet\pi_*{\cal O}}=\Lambda^\bullet \mtr E^*
\end{equation}
and the relative tangent bundle is given by
\begin{equation}
\mtr{T\pi}=\pi^*\mtr E^*\,\,.
\end{equation}
See also \cite[Th. VI.1.1]{FC}, where these properties are extended to toroidal
compactifications. The underlying real vector bundle of $E_\BC$ is
flat, as the period lattice determines a flat structure. Thus the topological
Pontrjagin classes
$p_j(E_\BC)$ vanish. For an action of $G=\mu_N$ on $X$  \cite[conjecture 3.2]{KR2}
states
\begin{conj}
$$
\ar\ch_G(\mtr{R^\bullet\pi_*{\cal O}})-a(T_g(\pi_\BC,\mtr{\cal O}))
=\pi^G_*\big(\ar\Td_G(\mtr{T\pi})(1-a(R_g(T\pi_\BC)))\big)
$$
\end{conj}
where $R_g$ denotes Bismut's equivariant $R$-class.
We shall assume henceforth that this conjecture holds.
Thus we obtain the equation
\begin{equation}\label{blah}
\ar\ch_G(\Lambda^\bullet \mtr E^*)-a(T_g(\pi_\BC,\mtr{\cal O}))
=\pi^G_*\big(\ar\Td_G(\pi^*\mtr E^*)(1-a(R_g(\pi^*E_\BC^*)))\big)\,\,.
\end{equation}
Using the equation
$$
\ar\ch_G(\Lambda^\bullet \mtr E^*)
=\frac{\ar c_{\rm top}(\mtr E^G)}{\ar \Td_G(\mtr E)}
$$
(\ref{blah}) simplifies to
$$
\frac{\ar c_{\rm top}(\mtr E^G)}{\ar \Td_G(\mtr E)}
-a(T_g(\pi_\BC,\mtr{\cal O}))
=\ar\Td_G(\mtr E^*)(1-a(R_g(E_\BC^*)))\pi^G_*\pi^*1
$$
or, using that $a(\ker \bar\partial\partial)$ is an ideal of square zero,
\begin{equation}\label{holla}
\ar c_{\rm top}(\mtr E^G)(1+a(R_g(E_\BC^*)))
-a(T_g(\pi_\BC,\mtr{\cal O})\Td_g(\mtr E_\BC))
=\ar\Td_G(\mtr E)\ar \Td_G(\mtr E^*)\pi^G_*\pi^*1\,\,.
\end{equation}
If $G$ acts fibrewise with isolated fixed points (over $\BC$), the left hand side of
equation (\ref{holla}) is an element of $\ar{CH}^0(Y)_{\BQ(\zeta_n)}+a({\rm
ker}\,\partial\mtr\partial)$. If $G$ does not act with isolated fixed points,
then the right hand side vanishes, $c_{\rm top}(E^G)$ vanishes and we
find \begin{equation}
\ar c_{\rm top}(\mtr E^G)
=a(T_g(\pi_\BC,\mtr{\cal O})\Td_g(\mtr E_\BC))
\,\,.
\end{equation}
As was mentioned in \cite[eq. (7.8)]{Ktori}, one finds in particular
$$
\ar c_d(\mtr E)=a(\gamma)\,\,.
$$
Now we restrict ourself to the action of the automorphism $(-1)$. We need to assume
that this automorphism corresponds to a $\mu_2$-action. This condition can always be
satisfied by changing the base Spec $D$ to Spec $D[\frac1{2}]$
(\cite[Introduction]{KR1} or \cite[section 2]{KR4}).
\begin{theor} \label{ful} Let $\pi:X\to Y$ denote
a principally polarized  abelian scheme of relative dimension $d$ over an arithmetic
variety $Y$. Set $\mtr
E:=(R^1\pi_* {\cal O},\|\cdot\|_{L^2})^*$. Assume \cite[conjecture 3.2]{KR2}. Then
the Pontrjagin classes of $\mtr E$ are given by
\begin{equation}
\ar p_k(\mtr
E)=(-1)^k\left(\frac{2\zeta'(1-2k)}{\zeta(1-2k)}+\sum_{j=1}^{2k-1}\frac1{j}
-\frac{2 \log 2}{1-4^{-k}}\right)(2k-1)!\,a(\ch(E)^{[2k-1]})\,\,.
\end{equation}
\end{theor}
The log 2-term actually vanishes in the arithmetic Chow ring over ${\rm
Spec}\,D[1/2]$.
\beginProof
Let $Q(z)$ denote
the power series in $z$ given by the Taylor expansion of
$$
4(1+e^{-z})^{-1}(1+e^z)^{-1}=\frac1{\cosh^2\frac{z}2}
$$
at $z=0$. Let $\ar Q$ denote the associated multiplicative arithmetic
characteristic class. Thus by definition for $G=\mu_2$
$$
4^d\ar\Td_G(\mtr E)\ar \Td_G(\mtr E^*)=\ar Q(\mtr E)
$$
and $\ar Q$ can be represented by Pontrjagin classes, as the power series $Q$ is
even. Now we can apply Lemma \ref{hups} for Pontrjagin classes to equation
(\ref{holla}). By a formula by Cauchy \cite[\S1, eq. (10)]{Hirz}, the summand of $\ar
Q$ consisting only of single Pontrjagin classes is given by taking the Taylor series
in
$z$ at $z=0$ of
$$
Q(\sqrt{-z})\frac{d}{dz}\frac{z}{Q(\sqrt{-z})}=
\frac{\frac{d}{dz}(z\cosh^2\frac{\sqrt{-z}}2)}{\cosh^2\frac{\sqrt{-z}}2}
=1+\frac{\sqrt{-z}}2\tanh\frac{\sqrt{-z}}2
$$
and replacing every power $z^j$ by $\ar p_j$. The bundle $\mtr
E$ is trivial, hence $\ar c_{\rm
top}(\mtr E^G)=1$. Thus we obtain by equation (\ref{holla})
with $\pi_*^G\pi^*1=4^d$
$$
\sum_{k=1}^\infty \frac{(4^k-1)(-1)^{k+1}}{(2k-1)!}\zeta(1-2k)\ar p_k
(\mtr
E)=-a(R_g(E_\BC))\,\,.
$$
Consider the zeta function $L(\a,s)=\sum_{k=1}^\infty k^{-s}\a^k$ for Re
$s>1$, $|\a|=1$. It has a meromorphic continuation to $s\in\BC$ which shall be
denoted by
$L$, too. Then
$L(-1,s)=(2^{1-s}-1)\zeta(s)$ and the function
$$
\widetilde R(\a,x):=\sum_{k=0}^\infty 
\left(\frac{\partial L}{\partial s}(\a,-k)+
L(\a,-k)\sum_{j=1}^k\frac{1}{2j}\right)\frac{x^k}{k!}
$$
defining the
Bismut equivariant $R$-class in \cite[Def. 3.6]{KR1} verifies for $\a=-1$
\begin{eqnarray}\nonumber
\widetilde R(-1,x)-\widetilde R(-1,-x)&=&
\sum_{k=1}^\infty \Big[(4^k-1)(2\zeta'(1-2k)+\zeta(1-2k)\sum_{j=1}^{2k-1}\frac1{j})
\\&&
-2\log 2\cdot 4^k\zeta(1-2k)\Big]\frac{x^{2k-1}}{(2k-1)!}
\,\,.
\end{eqnarray}
Thus we finally obtain the desired result.
\endProof
The first Pontrjagin classes are given by
$$
\ar p_1=-2 \ar c_2+\ar c_1^2,\qquad \ar p_2=2\ar c_4-2\ar c_3\ar c_1+\ar c_2^2,
\qquad \ar p_3=-2\ar c_6+2\ar c_5\ar c_1-2\ar c_4\ar c_2+\ar c_3^2\,\,.
$$
In general, $\ar p_k=(-1)^k 2 \ar c_{2k}+$products of Chern
classes. Thus knowing the Pontrjagin classes allows us to express the Chern classes
of even degree by the Chern classes of odd degree.
\section{A Hirzebruch proportionality principle and other applications}
Let $U$ denote the additive characteristic class associated to the power series
$$
\sum_{k=1}^\infty \left(\frac{\zeta'(1-2k)}{\zeta(1-2k)}+ \sum_{j=1}^{2k-1}\frac1{2j}
-\frac{\log 2}{1-4^{-k}}\right)
\frac{x^{2k-1}}{(2k-1)!}\,\,.
$$
\begin{cor} The part of $\ar\ch(\mtr E)$ in
$\ar{CH}^{\rm even}(Y)_\BQ$ is given by the formula
$$
\ar\ch(\mtr E)^{[{\rm even}]}=d-a(U(E))\,\,.
$$
\end{cor}
\beginProof
The part of $\ar\ch(\mtr E)$ of even degree equals
$$
\ar\ch(\mtr E)^{[{\rm even}]}=\frac1{2}\ar\ch(\mtr E\oplus\mtr E^*)\,\,,
$$
thus it can be expressed by Pontrjagin classes.
More precisely by Newton's formulae (\cite[\S 10.1]{Hirz}),
$$
(2k)!\ar\ch^{[2k]}-\ar p_1\cdot (2k-2)!\ar\ch^{[2k-2]}+\cdots
+(-1)^{k-1}\ar p_{k-1}2!\ar\ch^{[2]}
=(-1)^{k+1}k\ar p_k
$$
for $k\in\BN$.
As products of the arithmetic Pontrjagin classes vanish in
$\ar{CH}(Y)_\BQ$ by Lemma \ref{ful}, we thus observe that the part of $\ar\ch(\mtr
E)$ in
$\ar{CH}^{\rm even}(Y)_\BQ$ is given by
$$
\ar\ch(\mtr E)^{[{\rm even}]}=d+\sum_{k>0}\frac{(-1)^{k+1}\ar
p_k(\mtr E)}{2(2k-1)!}\,\,.
$$
Thus the result follows from Lemma \ref{ful}.
\endProof

\begin{theor}\label{height} Assume \cite[conjecture 3.2]{KR2}. 
There is a real number $r_d\in\BR$ and a Chern-Weil form $\phi(\mtr E)$ on
$Y_{\BC}$ of degree $(d-1)(d-2)/2$ such that
$$
\ar c_1^{1+d(d-1)/2}(\mtr E)=a(r_d\cdot c_1^{d(d-1)/2}(E)+\phi(\mtr E)\gamma)\,\,.
$$
\end{theor}
The form $\phi(\mtr E)$ is actually a polynomial with
integral coefficients in the Chern forms of $\mtr E$.
\beginProof Consider the graded ring $R_d$ given by $\BQ[u_1,\dots,u_d]$
divided by the relations
\begin{equation}\label{rela}
(1+\sum_{j=1}^{d-1}u_j)(1+\sum_{j=1}^{d-1}(-1)^ju_j)=1,\qquad u_d=0
\end{equation}
where $u_j$ shall have degree $j$ ($1\leq j\leq d$). This ring is finite
dimensional as a vector space over $\BQ$ with basis
$$
u_{j_1}\cdots u_{j_m},\qquad
1\leq j_1<\cdots<j_m<d\,\,, 1\leq m<d\,\,.
$$
In particular, any element of $R_d$ has degree $\leq\frac{d(d-1)}2$. As the
relation (\ref{rela})  is verified for $u_j=\ar c_j(\mtr E)$ up to multiples of
the Pontrjagin classes and $\ar c_d(\mtr E)$, any polynomial in the $\ar
c_j(\mtr E)$'s can be expressed in terms of the $
\ar p_j(\mtr E)$'s and $\ar c_d(\mtr E)$ if the corresponding polynomial in the
$u_j$'s vanishes in $R_d$.

Thus we can express $
\ar c_1^{1+d(d-1)/2}(\mtr E)$ as the image under $a$ of a topological
characteristic class of degree $\frac{d(d-1)}2$ plus $\gamma$ times a
Chern-Weil form of degree $\frac{(d-1)(d-2)}2$. As any element of
degree $\frac{d(d-1)}2$ in $R_d$ is proportional to $u_1^{d(d-1)/2}$, the
theorem follows.
\endProof
Any other arithmetic characteristic class of $\mtr E$ vanishing in $R_d$ can be
expressed in a similar way.

{\bf Example:} We shall compute $\ar c_1^{1+d(d-1)/2}(\mtr E)$ explicitly for
small $d$. 
Define topological cohomology classes $r_j$ by $\ar p_j(\mtr E)=a(r_j)$ via
Lemma \ref{ful}.
For $d=1$, clearly $$
\ar c_1(\mtr E)=a(\gamma)\,\,.
$$ 
In the case 
$d=2$ we find by the formula for $\ar p_1$
$$
\ar c_1^2(\mtr E)=a(r_1+2\gamma)=a\Big[(-1+\frac8{3}\log
2+24\zeta'(-1))c_1(E)+2\gamma\Big]\,\,.
$$
Combining the formulae for the first two Pontrjagin classes we get
$$
\ar p_2=2\ar c_4-2\ar c_3\ar c_1+\frac1{4}\ar c_1^4-\frac1{2} \ar c_1^2\ar
p_1+\frac1{4}\ar p_1^2.
$$
Thus for $d=3$ we find, using $c_3(E)=0$ and $c_1^2(E)=2c_2(E)$,
\begin{eqnarray}
\ar c_1^4(\mtr E)&=&a(2 c_1^2(E)r_1+4r_2+8c_1(E)\gamma)
\nonumber\\&=&\nonumber
a\Big[
(-\frac{17}{3}+\frac{48}{5}\log
2+48\zeta'(-1)-480\zeta'(-3))c_1^3(E)+8c_1(\mtr E)\gamma\Big]\,\,.
\end{eqnarray}
For $d=4$ one obtains
\begin{eqnarray*}
\ar c_1^7(\mtr E)&=&a\Big[ 64 c_2(E)c_3(E) r_1-(8 c_1(E)c_2(E)+32 c_3(E)) 
r_2+64 c_1(E) r_3\\
&&+ 16(7 c_1(\mtr E)c_2(\mtr E)-4 c_3(\mtr E))\gamma\Big]\,\,.
\end{eqnarray*}
As in this case $\ch(E)^{[1]}=c_1(E)$, $3!\ch(E)^{[3]}=-c_1^3(E)/2+3c_3(E)$ and
$5!\ch(E)^{[5]}=c_1^5(E)/16$, we find
\begin{eqnarray*}
\ar c_1^7(\mtr E)&=&a\Big[(-\frac{1063}{60}+\frac{1520}{63}\log
2+96\zeta'(-1)-600\zeta'(-3)+2016\zeta'(-5))c_1^5(E)
\\&&
+ 16(7 c_1(\mtr E)c_2(\mtr E)-4c_3(\mtr E))\gamma\Big]\,\,.
\end{eqnarray*}
\medskip

Now we are going to formulate an Arakelov version of Hirzebruch's
proportionality principle. In \cite[p. 773]{H2} it
is stated as follows: Let $G/K$ be a non-compact irreducible symmetric space
with compact dual $G'/K$ and let $\Gamma\subset G$ be a cocompact
subgroup such that $\Gamma\backslash G/K$ is a smooth manifold. Then there is an ring
monomorphism
$$
h:H^*(G'/K,\BQ) \to H^*(\Gamma\backslash G/K,\BQ)
$$
such that $h(c(TG'/K))=c(TG/K)$ (and similar for other bundles $F'$, $F$
corresponding to $K$-representation $V'$, $V$ dual to each other). This
implies in particular that Chern numbers on $G'/K$ and $\Gamma\backslash G/K$
are proportional\cite[p. 345]{H1}. Now in our case think about $Y$ as the
moduli space of principally polarized abelian varieties of dimension $d$. Its
projective dual is the Lagrangian Grassmannian $B_d$ over Spec $\BZ$
parametrizing isotropic subspaces in symplectic vector spaces of dimension
$2d$ over any field, $B_d(\BC)={\bf Sp}(d)/{\bf U}(d)$. But as the moduli space
is a non-compact quotient, the proportionality principle must be altered
slightly by considering Chow rings modulo certain ideals corresponding to 
boundary components in a suitable compactification. For that reason
we consider the Arakelov Chow group ${\rm CH}^*(\mtr B_{d-1})$, which is the
quotient of ${\rm CH}^*(\mtr B_d)$ modulo the ideal $(\ar c_d(\mtr S),
a(c_d(\mtr S)))$ with $\mtr S$ being the tautological bundle on $B_d$,
and we map it to $\ar{\rm CH}^*(Y)/(a(\gamma))$. Here $B_{d-1}$ shall be equipped
with the canonical symmetric metric. For the Hermitian symmetric space $B_{d-1}$, the
Arakelov Chow ring is a subring of the arithmetic Chow ring $\ar{\rm CH}(B_{d-1})$
(\cite[5.1.5]{GS}) such that the quotient abelian group depends only on
$B_{d-1}(\BC)$.
\begin{theor} Assume \cite[conjecture 3.2]{KR2}. 
There is a ring homomorphism
$$
h:{\rm CH}^*(\mtr B_{d-1})_\BQ\to \ar{\rm CH}^*(Y)/(a(\gamma))_\BQ
$$
with
$$
h(\ar c(\mtr S))=\ar c(\mtr E)
\left(1+a \left(\sum_{k=1}^{d-1}(\frac{\zeta'(1-2k)}{\zeta(1-2k)} -\frac{\log
2}{1-4^{-k}})(2k-1)!\ch^{[2k-1]}(E)\right)\right) $$
and
$$h(a(c(\mtr S)))=a(c(\mtr E))\,\,.$$
\end{theor}
Note that $S^*$ and $E$ are ample. One could as well map $a(c(\mtr S^*))$ to
$a(c(E))$, but the correction factor for the arithmetic characteristic classes
would have additional harmonic number terms.
\beginProof
The Arakelov Chow ring ${\rm
CH}^*(\mtr B_{d-1})$ has been investigated by Tamvakis in \cite{Tam}. Consider the
graded commutative ring $$
\BZ[\ar u_1,\dots,\ar u_{d-1}]\oplus \BR[u_1,\dots,u_{d-1}]
$$
where the ring structure is such that
$\BR[u_1,\dots,u_{d-1}]$ is an ideal of square zero. Let
$\ar R_d$ denote the quotient of this ring by the relations
$$
(1+\sum_{j=1}^{d-1}u_j)(1+\sum_{j=1}^{d-1}(-1)^ju_j)=1
$$
and
\begin{equation}\label{HarryRel}
(1+\sum_{k=1}^{d-1}\ar u_k)(1+\sum_{k=1}^{d-1}(-1)^k\ar u_k)=1-
\sum_{k=1}^{d-1}\left(\sum_{j=1}^{2k-1}\frac1{j}\right)(2k-1)!\ch^{[2k-1]}
(u_1,\dots,u_{d-1})
\end{equation}
where $\ch(u_1,\dots,u_{d-1})$ denotes the Chern character polynomial in the
Chern classes, taken of $u_1,\dots,u_{d-1}$.  Then by
\cite[Th. 1]{Tam}, there is a ring isomorphism $\Phi:\ar R_d\to {\rm CH}^*(\mtr
B_{d-1})$ with
$\Phi(\ar u_k)=\ar c_k(\mtr S^*)$, $\Phi(u_k)=a(c_k(\mtr S^*))$. The Chern
character term in (\ref{HarryRel}), which could be written more carefully as
$(0,\ch^{[2k-1]} (u_1,\dots,u_{d-1}))$, is thus mapped to $a(\ch^{[2k-1]}
(c_1(\mtr S^*),\dots,c_{d-1}(\mtr S^*)))$.  When writing the relation
(\ref{HarryRel}) as
$$
\ar c(\mtr S)\ar c(\mtr S^*)=1+a(\epsilon_1)
$$
and the relation in theorem \ref{ful} as
$$
\ar c(\mtr E)\ar c(\mtr E^*)=1+a(\epsilon_2)
$$

we see that a ring homomorphism $h$ is given by
$$
h(\ar c_k(\mtr S))=\sqrt{\frac{1+h(a(\epsilon_1))}{1+a(\epsilon_2)}}\ar c_k(\mtr E)
=(1+\frac1{2}h(a(\epsilon_1))-\frac1{2}a(\epsilon_2))\ar c_k(\mtr E)
$$
(where $h$ on im $a$ is defined as in the theorem).
Here the factor $1+\frac1{2}h(a(\epsilon_1))-\frac1{2}a(\epsilon_2)$ has even
degree, and thus
$$
h(\ar c_k(\mtr S^*))=\sqrt{\frac{1+h(a(\epsilon_1))}{1+a(\epsilon_2)}}\ar c_k(\mtr
E^*)
$$
which provides the compatibility with the cited relations.
\endProof
Note that this proof does not make use of the remarkable fact that
$h(a(\epsilon_1^{(k)}))$ and $a(\epsilon_2^{(k)})$ are proportional forms for any
degree
$k$.

In particular Tamvakis' height formula
\cite[Th. 3]{Tam} provides a combinatorial formula for the real number $r_d$
occurring in theorem
\ref{height}. Replace each term ${\cal H}_{2k-1}$ occurring in \cite[Th.
3]{Tam} by $$
-\frac{2\zeta'(1-2k)}{\zeta(1-2k)}-\sum_{j=1}^{2k-1}\frac1{j}
+\frac{2 \log 2}{1-4^{-k}}
$$
and divide the resulting value by half of the degree of $B_{d-1}$. Using
Hirzebruch's formula
$$
{\rm deg}\,B_{d-1}=\frac{(d(d-1)/2)!}{\prod_{k=1}^{d-1}(2k-1)!!}
$$
for the degree of $B_{d-1}$ (see \cite[p. 364]{H1}) and
the $\BZ_+$-valued function $g^{[a,b]_{d-1}}$ from \cite{Tam} counting
involved combinatorial diagrams, we obtain
\begin{cor}
The real number $r_d$ occurring in theorem \ref{height} is given by
\begin{eqnarray*}
r_d&=&\frac{2^{1+(d-1)(d-2)/2}\prod_{k=1}^{d-1}(2k-1)!!} {(d(d-1)/2)!}
\\&&\cdot
\sum_{k=0}^{d-2} 
\left(-\frac{2\zeta'(-2k-1)}{\zeta(-2k-1)}-\sum_{j=1}^{2k+1}\frac1{j}
+\frac{2 \log 2}{1-4^{-k-1}}
\right)
\\&&\cdot
\sum_{b=0}^{\min\{k,d-2-k\}}(-1)^b 2^{-\delta_{b,k}} g^{[k-b,b]_{d-1}}
\end{eqnarray*}
where $\delta_{b,k}$ is Kronecker's $\delta$.
\end{cor}
In \cite[Th. 2.5]{Geer} van der Geer shows that $R_d$ embeds into the
(classical) Chow ring ${\rm CH}^*({\cal A}_d)_\BQ$ of the moduli stack ${\cal A}_d$
of principally polarized abelian varieties. Using this result one finds

\begin{lemma}\label{Harry}
Let $Y$ be a regular finite covering of the moduli space ${\cal A}_d$ of principally
polarized abelian varieties of dimension $d$. Then for any non-vanishing polynomial
$p(u_1,\dots, u_{d-1})$ in $R_d$,
$$
h(p(\ar c_1(\mtr S),\dots ,\ar c_{d-1}(\mtr S)))\notin {\rm im}\,a\,\,.
$$
In particular, $h$ is non-trivial in all degrees. Furthermore, $h$ is injective iff
$a(c_1(E)^{d(d-1)/2})\neq0$ in $\ar{\rm CH}^{d(d-1)/2+1}(Y)_\BQ/(a(\gamma))$.
\end{lemma}
\beginProof
Consider the canonical map $\zeta:\ar{\rm CH}^*(Y)_\BQ/(a(\gamma))\to {\rm
CH}^*(Y)_\BQ$. Then
$$\zeta(h(p(\ar c_1(\mtr S), \dots ,\ar c_{d-1}(\mtr S)))) = p( c_1(E), \dots
,c_{d-1}(E))\,\,,$$ and the latter is non-vanishing according to
\cite[Th. 1.5]{Geer}. This shows the first part.

If $a(c_1(E)^{d(d-1)/2})\neq0$ in $\ar{\rm CH}^{d(d-1)/2+1}(Y)_\BQ/(a(\gamma))$,
then by the same induction argument as in the proof of \cite[Th. 2.5]{Geer} $R_d$
embeds in
$a(\ker \bar\partial\partial)$. Finally, by \cite[th. 2]{Tam} any element $z$ of
$\ar R_d$ can be written in a unique way as a linear combination of
$$
\ar u_{j_1}\cdots \ar u_{j_m}\quad{\rm and}\quad
u_{j_1}\cdots u_{j_m},\qquad
1\leq j_1<\cdots<j_m<d\,\,, 1\leq m<d\,\,.
$$
Thus if $z\notin{\rm im}\,a$, then $h(z)\neq0$ follows by van der Geer's
result, and if $z\in{\rm im}\,a\setminus\{0\}$, then $h(z)\neq0$ follows by embedding
$R_d\otimes
\BR$.
\endProof

\end{document}